\documentclass[10pt]{amsart}
\usepackage{amsmath, latexsym, 
}
\usepackage{amssymb,amsthm}
\usepackage{amscd,graphics}

\newtheorem{proposition}{Proposition}[section]

\newtheorem{theorem}[proposition]{Theorem}

\theoremstyle{remark}

\theoremstyle{definition}
\newtheorem*{definition}{Definition}
\def\real{\mathbb{R}}

\def\complex{\mathbb{C}}

\def\var{\mathrm{var}}

\def\LL{\mathcal{L}}

\begin{document}
\title[Analyticity of the SRB measure for
families of Collet-Eckmann maps]{Analyticity of the SRB measure for holomorphic families
of quadratic-like Collet-Eckmann maps}
\author{Viviane Baladi and Daniel Smania} 
\address{D.M.A., UMR 8553,\'Ecole Normale Sup\'erieure,  75005 Paris, France}
\email{viviane.baladi@ens.fr}
\begin{abstract}
We show that if $f_t$ is a holomorphic family of quadratic-like maps with all
periodic orbits repelling so that for each real $t$ the map $f_t$ is a
real Collet-Eckmann $S$-unimodal map then, writing
$\mu_t$ for the unique absolutely continuous invariant probability measure
of  $f_t$, the map
$$
t\mapsto \int \psi \, d\mu_t
$$
is real analytic for any real analytic function $\psi$.
\end{abstract}
\address{
Departamento de Matem\'atica,
ICMC-USP, Caixa Postal 668,  S\~ao Carlos-SP,
CEP 13560-970
S\~ao Carlos-SP, Brazil}
\email{smania@icmc.usp.br}
\date{\today} 
\thanks{V.B. is partially supported by ANR-05-JCJC-0107-01.
D.S. is partially supported by CNPq 470957/2006-9
and 310964/2006-7, FAPESP 2003/03107-9. 
D. S. thanks the DMA of \'Ecole Normale Sup\'erieure for
hospitality during a visit
where a crucial part of this work was done.
V.B. wrote part of this paper while visiting the Universidad Cat\'olica del
Norte, Antofagasta, Chile, whose hospitality is gratefully acnowledged.
We thank D. Sands for very helpful comments.}
\maketitle

%%%%%%%%%%%%%%%%%%%%%

\section{Introduction and statement of the theorem}

If $t\mapsto f_t$ is a smooth one-parameter family of dynamics 
$f_t$ so that $f_0$ admits
a unique SRB measure $\mu_0$, it is natural to ask whether the map $t\mapsto \mu_t$,
where $t$ ranges over a set $\Lambda$
of parameters such that $f_t$ has (at least) one SRB measure 
$\mu_t$,  is differentiable at $0$ (in the sense of Whitney if $\Lambda$ does
not contain a neighbourhood of $0$, as suggested by Ruelle \cite{Rupr}).
Katok, Knieper, Pollicott, and Weiss \cite{KKPW} gave a positive answer to
this question in the setting of $C^3$ families
of transitive Anosov flows (here, $\Lambda$ 
is a neighbourhood of $0$), 
showing that $t\mapsto \int \psi\, d\mu_t$ is differentiable,
for all smooth $\psi$.
If $f_0$ is a $C^3$ mixing Axiom A attractor and the family $t\mapsto f_t$
is $C^3$, Ruelle \cite{Ruok} not only proved that $t \mapsto \int \psi\, d\mu_t$ is
differentiable, but also gave an explicit formula (the linear response formula)
for the derivative. Ruelle \cite{Rupr} suggested that this formula, appropriately
interpreted, should hold in much greater generality.
Indeed, Dolgopyat \cite{Dolgo} obtained the linear response formula for
a class of partially hyperbolic diffeomorphisms.
In a previous work \cite{BS1, BS2}, we found that in the (non structurally
stable) setting of
piecewise expanding unimodal interval maps, the SRB measure is differentiable
if and only if the path $f_t$ is tangent to the topological class of $f_0$,
that is, if and only if $\partial_t f_t|_{t=0}$ is horizontal.
When differentiability holds, Ruelle's candidate for the derivative,
as interpreted in \cite{Ba}, gives the linear response formula.
 (We refer to \cite{Ba,BS1, BS2},
which also contain conjectures about smooth, not necessarily analytic,
Collet--Eckmann maps, for more information and additional references.)
Then, Ruelle \cite{RuStr} proved the linear response formula for a class
of nonrecurrent \footnote{I.e., $\inf_k d(f^k_t(c),c) >0$,
where $c$ denotes the critical point.} analytic unimodal interval maps $f_t$,
assuming that all $f_t$ stay in the topological class of $f_0$.
In the present work, we consider holomorphic (that is, complex analytic)
families $f_t$ of quadratic-like holomorphic Collet--Eckmann maps.
Our assumptions imply (using classical holomorphic motions)
that all $f_t$ lie in the same conjugacy class. Generalising one of the
arguments in \cite{KKPW}, we are able to show that 
$
t\mapsto \int \psi \, d\mu_t
$
is real analytic for any real analytic function $\psi$.

Let us now state our result more precisely.
Let $I=[-1,1]$.
A $C^3$ map $f:I \to I$ is an
{\it $S$-unimodal} map if it has $c=0$ as unique critical point,
and $f$ has nonpositive Schwarzian derivative, that is
$\frac{f'''}{f'}- \frac{3}{2}\bigl ( \frac{f''}{f'}\bigr )^2\le 0$
except at $c$.
An $S$-unimodal map is called Collet-Eckmann  if there 
exist $C >0$ and $\lambda_c >1$ so that 
$|(f^n)'(f(c))|\ge C \lambda_c^n$ for all $n \ge 1$.
In this paper, we shall only consider $S$-unimodal maps with
$f''(c)\ne 0$.

In Section~\ref{proo} we shall define precisely the notion of
a
holomorphic  (complex analytic) family of {\it quadratic-like} maps
{\it in a neighbourhood of $I$} and prove the main result of this work:

\begin{theorem} 
Let $t \mapsto f_t$ be a holomorphic family of quadratic-like
maps in a neighbourhood of $I$,
with all periodic orbits repelling. Assume in addition that  for each small real $t$ the map
$f_t$  restricted to  $I$ is a (real) Collet-Eckmann S-unimodal map.
Then there exists $\epsilon >0$
so that for each  real analytic  $\psi:I\to \complex$,
the map 
$$t\mapsto \int \psi \rho_t \, dx,
$$
where $\rho_t$ is the invariant
density of $f_t$, is  real analytic on $(-\epsilon,\epsilon)$.
\end{theorem}

The quadratic-like assumption implies that $f''_t(c)<0$.
The fact that periodic orbits are repelling implies that $f_t$ is topologically conjugated
with $f_0$ (see our use of  Ma\~n\'e-Sad-Sullivan \cite{MSS}
in the beginning of the proof of the theorem in Section~\ref{proo}).
Besides Ma\~n\'e-Sad-Sullivan \cite{MSS} the other main
ingredient of our proof are the results and constructions
of Keller and Nowicki \cite{KN}\,  which allow us
to exploit dynamical zeta functions, following the argument
in the work of Katok--Knieper--Pollicott--Weiss \cite[First proof of Theorem ~ 1]{KKPW}.

The extension from quadratic-like to polynomial-like is
straightforward, and we stick to the nondegenerate case
$f''(c)\ne 0$  for the sake of simplicity of exposition.
As the proof uses only real-analyticity of the holomorphic motions
$t\mapsto h_t$, it is conceivable that the conclusion of the theorem holds
if $f_t $ is a real analytic family of quadratic-like maps, using
ideas of \cite{ALM}, but this generalisation  appears to
be nontrivial.

\section{Proof of the Theorem}
\label{proo}

Before we prove the theorem, let us define precisely the objects we are studying:

\begin{definition}
We say that $f_t$ is a {\it holomorphic family
of quadratic-like maps in a neighbourhood
of $I$} if there exists
a complex neigbourhood $U$  of $I$
so that $t\mapsto f_t$  is a holomorphic map
from a complex neighbourhood of zero  to the Banach space $B(U)$ of
holomorphic functions on $U$ extending continously  to  
$\overline U$ (with the supremum norm), such that:
\begin{itemize}
\item
For real $t$, the map $f_t$ is real on $\Re U$, with
$f_t(I)\subset  I$ and $f_t(-1)=f_t(1)=-1$.
\item
There exist simply connected complex domains  $W$ and $V$, 
whose boundaries are analytic Jordan curves,
with
$
I\subset W$, $I\subset V$, $\overline V \subset U$,
$\overline V \subset W$,
and  so that  $f_0: V\mapsto W$
is  a double-branched ramified covering, with $c=0$ as a unique 
critical point. (That is, $f_0:V\mapsto W$ is a  quadratic-like restriction  
of $f_0$.)
\end{itemize}
\end{definition}

If $f_t$ is a holomorphic family
of quadratic-like maps in a neighbourhood
of $I$  then it is easy to see 
\footnote{Indeed, $\partial W$ is an analytic Jordan curve, and $f_0$ has no critical 
point on $\partial  V$. If $f_t \in B(U)$ is close to $f_0$,
there is a simply connected domain $V_t$ close to $V$ such that 
$f_t(V_t)=W$, and the boundary of $\partial V_t$ is  a Jordan curve, by
the implicit function theorem. Then $f_t:V_t \to W$ is a 
quadratic-like extension.}
that for small complex $t$,  denoting
by $V_t$ the connected component of $f^{-1}_t (W)$ containing $0$, then 
$f_t: V_t \mapsto W$  is a quadratic-like restriction of $f_t$.
We may then give another definition:
 
\begin{definition}
We say that  $f_t$ is a holomorphic family
of quadratic-like maps in a  neighbourhood
of $I$ {\it with all periodic orbits repelling,} if
$f_t$ is a holomorphic family
of quadratic-like maps in a  neighbourhood
of $I$ so that, for each small complex $t$, the map
$f_t$ only has repelling periodic orbits in $V_t$.
\end{definition}

\begin{proof}
Since we assumed that
all periodic points of $f_t$ are repelling,  \cite[Theorem B]{MSS}
(the result there is quoted for polynomial maps, but the proof
immediately extends to polynomial-like) implies that
there exists  a
holomorphic motion of  the Julia set $K(f_0)$ of $f_0$, that is, a map
$h: D \times K(f_0) \to C$
where $D=\{ z \in \complex \mid |z|< \epsilon_0\}$
for some $\epsilon_0>0$, such that
 for each $x \in K(f_0)$ the map
$t\mapsto h_t(x)$
is holomorphic, and for every $t\in D$ the function
$x\mapsto h_t(x)$
is continuous and injective on $K(f_0)$, with
 $$h_t \circ f_0 = f_t \circ h_t\, .
 $$
 (In particular, $h_t$ is a homeomorphism from $K(f_0)$ to
 $K(f_t)$.)
Our assumptions imply that $[f^2_0(0), f_0(0)] =K(f_0)\cap \real$ and
$h_t(K(f_0)\cap \real)= K(f_t)\cap \real=[f^2_t(0), f_t(0)]$.
From now on, we only use real analyticity of $t \mapsto f_t(x)$ and
$t\mapsto h_t(x)$ for $x \in [f^2(0), f(0)]$.

We next claim that our assumptions guarantee that each $f_t$ satisfies the
technical requirement needed by Keller and Nowicki \cite[(1.2)]{KN}.
Denoting by $\var_J \phi$ the total variation of a function $\phi$
on an interval $J$, and writing $f=f_t$, we  need to check that there is
that  a constant $M>0$ such that:
\begin{enumerate}
\item\label{one}
$M^{-1} < \sup_I  \frac{|x-c|}{|f'(x)|} 
+ \var_I \frac{|x-c|}{|f'(x)|}   < M$,
\item\label{two}
$\var_{J_u}  \frac{|f(x)-f(u)|}{|x-u||f'(x)|}  < M$
where $J_u=[-1,u]$ if $u<c$ and $=[u,1]$ if $u >c$.
\end{enumerate}
Let $\delta_1>0$ be so that $|f''(y)|> |f''(c)|/2$ if $|y-c|<\delta_1$.
It suffices to prove (\ref{one}) and (\ref{two}) for $|x-c|<\delta_1$
and $|u-c|<\delta_1$, and we restrict to such
points. Noting that for every such $x\ne c$ there
exist  $y_x$, $z_x$, and $\tilde z_x$, between $x$ and $c$,
so that
$$
\frac{|x-c|}{|f'(x)|} =-\frac{x-c}{f'(x)-f'(c)} =
-\frac{1}{f''(y_x)}\, ,
$$
and (use $f''(x)= f''(c)+ f^{(3)}(z_x)(x-c)$ and $f'(x)=f''(c) (x-c)+ f^{(3)}(\tilde z_x) 
\frac{(x-c)^2}{2}$)
$$
\partial_x \frac{|x-c|}{|f'(x)|} =\frac{-f'(x)+(x-c)f''(x)}{(f'(x))^2} =
\frac{(x-c)^2}{(f'(x))^2} \bigl (f^{(3)}(z_x)-\frac{f^{(3)}(\tilde z_x)}{2}\bigr )\, ,
$$
the first two conditions hold because $f$ is $C^3$.
For the third condition, consider 
$x\ge u>c$ (the other case is symmetric).
Since 
$$
\frac{f(x)-f(u)}{(x-u)f'(x)}=
1 + \frac{x-u}{f'(x)} \frac{f''(z_x)}{2}=1 + \frac{x-u}{f'(x)} \frac{f''(z_x)}{2 f''(y_x)}
\, , 
$$
and 
$
0<-\frac{x-u}{f'(x)}<  -\frac{x-c}{f'(x)}
$,
we get that $\bigl | \frac{f(x)-f(u)}{(x-u)f'(x)}\bigr | $ is bounded on $[u,1]$, uniformly
in $u$. Finally, since
$$
\partial_x \frac{x-u}{f'(x)}  =\frac{f'(x)-(x-u)f''(x)}{(f'(x))^2} \, , 
$$
analyticity of $f$ implies that  $\partial_x \frac{x-u}{f'(x)}$ changes signs
finitely many times, uniformly in $u$, proving (\ref{two}).

Also, the results of Nowicki--Sands \cite{NS} and Nowicki--Przytycki \cite{NP}
ensure (see Appendix~\ref{Sands}) that there exist $\lambda_c >1$,
$\lambda_{per} >1$, $\lambda_\eta >1$, and $\epsilon_1 >0$ so that, for 
each $|t|< \epsilon_1$, there is $C_t >0$ with
\begin{equation}\label{unifCE}
|(f^n_t)'(f_t(0)|\ge C_t \lambda_c ^n \, , \forall n \ge 1\, ,
\end{equation}
and so that for each $x\in I$ so that $f_t^p(x)=x$ for some $p\ge 1$, we have
\begin{equation}\label{unif2}
|(f^p_t)'(x)|\ge C_t \lambda_{per}^p \, ,
\end{equation}
and, finally, setting
$$
\lambda_\eta(t):=\liminf_{n \to \infty}
\{ |\eta|^{-1/n} \mid \eta \subset I \mbox{ is the biggest monotonicity interval of } f_t^n\}\, ,
$$
\begin{equation}\label{unif3}
\inf_{|t|<\epsilon_1} \lambda_\eta(t) > \lambda_\eta\, .
\end{equation}
In other words, the hyperbolicity constants are uniform in
$t$, 
guaranteeing uniformity when applying
the results of Keller and Nowicki \cite{KN}.
(We choose $\epsilon_1 <\epsilon_0$.)

We now adapt the strategy used in the first proof of
\cite[Theorem 1]{KKPW}.  Fix $\psi$ and, for $x \in I$ so that
$f^p_0(x)=x$ for $p\ge1$, and for small real $s$ and
$t$, consider 
\begin{equation}\label{gst}
g_{s,t}(x)=  \frac{e^{ s \psi( h_t(x))}}{|f'_t ( h_t(x))|}\, .
\end{equation}
Since $\psi$ is real analytic, the analyticity of $t \mapsto h_t$
and of $t \mapsto f_t$ together with (\ref{unif2}) imply that there is
$\epsilon_2>0$ so that, for every periodic point $x\in I$ 
of period $p\ge 1$ for
$f$, the function
$$
(t,s) \mapsto  g_{s,t}^{(p)}(x):=
 \frac{e^{ s \sum_{k=0}^{p-1}\psi (h_t(f^k(x))}}{|(f^p_t)' ( h_t(x))|}
$$
is real analytic in $|s|<\epsilon_2$ and $|t|<\epsilon_2$, uniformly in $x$.
We take $\epsilon_2<\epsilon_1$.

Therefore, the dynamical zeta function defined by
\begin{align}
\zeta(s,t,z)&:=\exp \sum_{p = 1}^\infty \frac{z^p}{p}
\sum_{x \in I : f^p_0(x)=x}  g^{(p)}_{s,t}(x)
\end{align}
has the following property: There exists $\delta_2>0$ so that
for each  $|z|< \delta_2$ the function
$\zeta(s,t,z)$ is real analytic in $|t|< \epsilon_2$, $|s|< \epsilon_2$, and  
so that for each  $(s,t)$ with $|t|< \epsilon_2$, $|s|< \epsilon_2$
the map $\zeta(s,t,z)$ is holomorphic and nonvanishing in $|z|< \delta_2$.

Now, $h_t \circ f_0 = f_t \circ h_t$ immediately implies
\begin{equation}\label{zeta2}
\zeta(s,t,z)=\exp \sum_{p = 1}^\infty \frac{z^p}{p}
\sum_{y \in I : f^p_t(y)=y} \frac{e^{ s \sum_{k=0}^{p-1}\psi( f_t^k (y))}}
{|(f^p_t)' ( y)|}\, .
\end{equation}

Recall (\ref{unifCE}, \ref{unif2}, \ref{unif3})
and take $\Theta \in (0,1)$ with
$$
\Theta^{-1} <  \min \{\lambda_\eta, 
\sqrt {\min (\lambda_c, \lambda_{per}) }\}\, .
$$
Keller and Nowicki \cite[Theorem 2.1]{KN} prove that, if $\epsilon_3\in(0,\epsilon_2)$ is
small enough, then  for  $|s|<\epsilon_3$ and $|t|< \epsilon_3$
the transfer operator
\footnote{Our parameter $s$ is called $t$ in \cite{KN}, the parameter
$\beta$ in \cite{KN} is  $\beta=1$, and  our parameter
$t$ corresponds to changing the dynamics.}
$$
 \LL_{s,t} \varphi (x)=\sum_{\hat f_t(y)=x} \frac{ \omega_t(y)}
{\omega_t(x)} \frac{\exp( s \psi(y))}{|\hat f'_t ( y)|} \varphi(y)\, ,
$$
acting on functions of bounded variation on a suitable Hofbauer tower
extension $\hat f_t :\hat I
\to \hat I$ of $f_t$
\cite[Section 3]{KN},
endowed with an appropriate
\cite[\S 6.2]{KN} cocycle $\omega_t$ (which embodies the singularities
along the postcritical orbit of $f_t$), is a bounded operator. 
If $s=0$ then the spectral radius $\lambda_{0,t}$ of $\LL_{s,t}$
is equal to $1$,
it is a simple eigenvalue (whose eigenvector gives the invariant
density $\rho_t$ of $f_t$), and the rest of the spectrum is contained in
a disc of strictly smaller radius.
In addition, the essential spectral radius $\theta_{s,t}$
of  $\LL_{s,t}$ satisfies
$\sup_{|t|< \epsilon_3, |s| < \epsilon_3} \theta_{s,t} < \Theta$,
and  for each $|t|<\epsilon_3$ the spectral radius
\footnote{Note that $\lambda_{s,t}$ is the exponential of the topological pressure 
of $s\psi-\log |f'_t|$ for $f_t$, and that $\rho_t\, dx$
is the equilibrium state for $f_t$ and $-\log |f'_t|$.} 
$\lambda_{s,t}>\Theta$
of $\LL_{s,t}$ is an analytic function
\cite[Prop. 4.2]{KN} of $s$.
Also, perturbation theory gives (see \cite[(5.2)]{KN})
\begin{equation}\label{perturb}
\partial_s \log \lambda_{s,t} |_{s=0}=  \int \psi  \rho_t\, dx\, .
\end{equation}

Keller and Nowicki also
show \cite[Theorem 2.2]{KN}
that for $|t|< \epsilon_3$ and $|s|< \epsilon_3$ the power series
$\zeta(s,t,z)$ defined by (\ref{zeta2})
extends meromorphically to the disc of radius $\Theta^{-1}$
(where it does not vanish, by \cite[Prop. 4.3 and Lemma 4.5]{KN}), and its
poles $z_k$ in this disc are in bijection with the eigenvalues $\lambda_k$
of $\LL_{s,t}$, via $\lambda_k=z_k^{-1}$.
(The order of the zero coincides with the algebraic multiplicity
of the eigenvalue.)
It follows that $z \mapsto \zeta(s,t,z)^{-1}$ is holomorphic 
in the disc of radius $\Theta^{-1}$. This disc contains $\lambda_{s,t}^{-1}$, which
is a simple zero.

To end the proof, recalling (\ref{perturb}), it suffices to see that
 $(s,t)\mapsto \lambda_{s,t}$ is real analytic, but this easily
follows from Shiffman's \cite{Sh} real analytic Hartogs' theorem
(see Appendix~B or \cite[Thm p. 589]{KKPW}) applied to $d(s,t,z)= \zeta(s,t,z)^{-1}$, which implies that for each $(s,t)\in  (-\epsilon_3, \epsilon_3)\times
(-\epsilon_3, \epsilon_3)$ the map $z\mapsto d(s,t,z)$ is holomorphic
in $|z|<\Theta^{-1}$.
Indeed, by the implicit function theorem,
the simple zeroes of  $d(s,t,\cdot)$ depend 
real analytically on $s$ and $t$.
(We used the same $\epsilon_i$ discs for the $s$ and $t$ variable,
but a more careful analysis shows that $\epsilon$ in the statement of the theorem
may be selected independently of $\psi$.)
\end{proof}

\appendix
\section{Uniformity of the hyperbolicity constants}
\label{Sands}

We start with a preliminary observation\footnote{We thank Duncan Sands for his explanations.}:
Let $g$ be an $S$-unimodal Collet--Eckman map (with $g''(0)< 0$, say).
Denote by $\lambda_c(g)$, $\lambda_{per}(g)$, and
$\lambda_\eta(g)$  the constants defined by (\ref{unifCE}, \ref{unif2}, \ref{unif3})
(replacing $f_t$ by $g$).
Nowicki and Sands \cite{NS}  proved that
if  $g$ is an $S$-unimodal map and $\lambda_{per}(g)>1$ then $\lambda_{c}(g)>1$.
A careful study of their proof shows that 
 $\lambda_c(g) > \lambda_{per}(g)^\alpha$,
where the exponent $\alpha >0$ only depends
on the maximum length $N(g)$ of ``almost-parabolic funnels" of $g$
(see \cite[Lemma~6.6]{NS} for a definition of $N(g)$, which
can be bounded by a function of $1/\log(\lambda_{per}(g))$
and $\sup |g'|$). Since $N(g)$ 
is in fact  invariant under topological conjugacy and
$f_t$ is topologically conjugated to $f_0$, we conclude that
$\lambda_c(f_t) > \lambda_{per}(f_t)^\alpha$, with
$\alpha >0$ uniform in small $t$.

Next, recall that Nowicki and Przytycki \cite{NP} proved that if $g$ and
$\tilde g$ are $S$-unimodal  maps (with $g''(c)\ne 0$
and $\tilde g''(c)\ne 0$, say)  conjugated by a homeomorphism
of the interval and  $g$ is Collet--Eckmann, then $\tilde g$ is Collet--Eckmann.
Take $g=f_0$ and $\tilde g=f_t$ (in particular, $f_t$ is $C^2$
close to $f_0$ and $t\mapsto h_t$ is smooth).
Then it is not very difficult to see that 
the constants $M=M(f_t)>0$, $P_4=P_4(f_t)>0$,
and $\delta_4=\delta_4(f_t)>0$ from the
topological characterisation  (``finite criticality")  of Collet--Eckmann in 
\cite[(4) p. 35]{NP})
are uniform in small $t$.

Recall that our assumptions imply $f_t''(c)\ne 0$ for all small $t$, so that
the constant denoted $l_c$ in \cite{NP} is $l_c=2$.
Section ~2 of \cite{NP}, and in particular the use of the Koebe principle there,
implies that there exists a (universal)
function 
$q:\real^+_*\times (0,1) \to (0,1)$
with
$q(M,1/4)<1/2$ for any   $M$ (see \cite[Lemma~2.2]{NP}),
and so that
$\lambda_{per}(f_t) > \bigl (1- 2q(M(f_t),1/4) \bigr )^{-1}$.
Therefore, 
$\lambda_{per}(f_t)>1$ is uniformly bounded away from $1$ for small $t$.
The preliminary observation then implies that
$\lambda_c(f_t)$ is also uniformly bounded in $t$.
By \cite[Proposition 3.2]{No}
(see also \cite[p. 35]{NP}), this implies a uniform lower bound for
$\lambda_\eta(f_t)$.
(Indeed, in the notations of \cite[\S3]{No}, we have
$\lambda_\eta=\lambda_5=\lambda_4\ge \lambda_3 =\lambda_1 \ge \sqrt{\lambda_c}$.)

\section{Shiffman's real analytic Hartogs' extension theorem}

\begin{theorem}\cite{Sh}
Let $\delta >0$ and $0 <r <R$.
Assume that 
$$d : (-\delta, \delta)^2 \times \{ z \in \complex \mid |z| <R\} \to \complex$$
satisfies the following conditions:
\begin{itemize}
\item
For each $(s,t)\in  (-\delta, \delta)^2$ the map $z\mapsto d(s,t,z)$ is holomorphic
in $|z|<R$.
\item
For each $|z| <r$ the map $(s,t)\mapsto d(s,t,z)$ is real
analytic in  $(-\delta, \delta)^2$.
\end{itemize}
Then $d(s,t,z)$ is real analytic on  $(-\delta, \delta)^2 \times \{  |z| <R\}$.
\end{theorem}

Note that the above theorem fails if real analyticity is replaced by
$C^k$ for $k\le \infty$.

The theorem holds because $|z| < r$ is not pluripolar in $|z|<R$.
Shiffman's result
is based on deep work of Siciak \cite{Si}

\bibliographystyle{amsplain}

\end{document}